\documentclass{amsart}
\newtheorem{theorem}{Theorem}[section]
\newtheorem{lemma}[theorem]{Lemma}

\theoremstyle{definition}

\newtheorem{remark}[theorem]{Remark}

\newtheorem{example}[theorem]{Example}

\newtheorem{cor}[theorem]{Corollary}
\theoremstyle{remark}

\begin{document}
\title{On  rectangular constant in normed linear spaces}

\author{Kallol Paul, Puja Ghosh  and Debmalya Sain}
\address{Department of Mathematics, Jadavpur University, Kolkata 700032, India}
\email{kalloldada@yahoo.co.in, saindebmalya@gmail.com}
 \thanks{The second and third author would like to thank CSIR and UGC respectively for the financial support.}

\subjclass[2010]{Primary 46B20, Secondary 47A30}
\keywords{ Birkhoff-James Orthogonality, rectangular constant.}
 
\date{}
\begin{abstract}
We study the properties of rectangular constant  $ \mu(\mathbb{X}) $  in a normed linear space $\mathbb{X}$. 
We prove that $ \mu(\mathbb{X}) = 3$ iff the unit sphere contains a straight line segment of length 2. In fact, we prove that the rectangular modulus attains its upper bound iff the unit sphere contains a straight line segment of length 2. 
 We  prove that if the dimension of the space $\mathbb{X}$ is finite then $\mu(\mathbb{X})$ is attained.  We also prove that a normed linear space is an inner product space iff  we have sup$\{\frac{1+|t|}{\|y+tx\|}$: $x,y \in S_{\mathbb{X}}$ with $x\bot_By\} \leq \sqrt{2}$ $\forall t$ satisfying $|t|\in (3-2\sqrt{2},\sqrt{2}+1)$.
\end{abstract}
\maketitle
\section{Introduction}

\noindent Let $(\mathbb{X}, \|~~\|)$ be a  normed linear space, unless otherwise mentioned we assume the space to be real. Let $ B_{\mathbb{X}}=\{x\in \mathbb{X} : \|x\|\leq 1\} $ and $ S_{\mathbb{X}}=\{ x \in \mathbb{X}: \|x\|=1 \} $ be the unit ball and unit sphere respectively of the normed linear space $\mathbb{X}$.  There are various notions of orthogonality in a normed linear space. Among them one of the most important and widely used notion is orthogonality in the sense of Birkhoff-James\cite{2,7}:\\
For any two elements $ x , y \in \mathbb{X}$, $x$ is said to be orthogonal to $y$ in the sense of Birkhoff-James, written as $x\bot_B y$ iff $\|x\| \leq \|x + \lambda y\|$ for all $\lambda \in K(= \mathbb{R} ~or~ \mathbb{C}) $.\\
In the year of $1969$, Joly \cite{8} introduced the notion of rectangular constant of a real normed linear space of dimension $\geq 2$, which is defined as follows: 
\[\mu(\mathbb{X})=\sup_{x\perp_B y}~\Big\{\frac{\|x\|+\|y\|}{\|x+y\|}\Big\}.\]  
The notion of rectangular constant plays a very important role in the  geometry of normed linear spaces and it has been studied by Joly\cite{8}, del Rio and Benitez\cite{3}, Desbiens\cite{4} and  Baronti and  Casini \cite{1}. In 1999 the concept of the rectangular constant of a normed linear space was generalized by  Serb \cite{9} to introduce  rectangular modulus, $ \mu_{\mathbb{X}}(\lambda)$, as a function 
$ \mu_{\mathbb{X}}:(0,\infty)\longrightarrow \mathbb{R}$ defined as follows 
\[ \mu_{\mathbb{X}}(\lambda) = \sup \Big\{\max\{\frac{\lambda^2+t}{\|\lambda u+tv\|},\frac{1+\lambda^2 t}{\|u+\lambda tv\|}\}: t> 0, u,v\in S_{\mathbb{X}}, u \perp_B v \Big\} .\]
 For $ \lambda = 1$,  $\mu_{\mathbb{X}}(1)$ is the rectangular constant of the normed linear space $\mathbb{X}$. He also defined $*$-rectangular modulus $ \mu^{*}_{\mathbb{X}}(\lambda) $ as 
\[ \mu^{*}_{\mathbb{X}}(\lambda) = \sup \Big\{\frac{\lambda^2+t}{\|\lambda u+tv\|}\}: t> 0, u,v\in S_{\mathbb{X}}, u \perp_B v\Big \} \]
and proved that 
\[ \mu_{\mathbb{X}}(\lambda)  = \max\{\mu^{*}_{\mathbb{X}}(\lambda), \lambda\mu^{*}_{\mathbb{X}}(1/\lambda)\}.\]
Serb \cite{9} proved that 
\[ \sqrt{1 + \lambda^2} \leq \mu_{\mathbb{X}}(\lambda) \leq \max\{ \lambda +2, 1 + 2 \lambda\} ~~ \forall \lambda > 0.\]
\noindent It is well known that for any real normed linear space $\mathbb{X},$ $\sqrt{2} \leq\mu(\mathbb{X})\leq 3.$  Joly\cite{8} proved that if dimension of $\mathbb{X}$ is greater than $2$ then   $\mu(\mathbb{X})=\sqrt{2}$ iff $\mathbb{X}$ is an inner product space. Later  del Rio and Benitez\cite{3} proved it for dimension 2.  However there is no characterisation for the attainment of upper bound of the rectangular constant. In this paper we prove that $ \mu(\mathbb{X}) = 3$ iff the unit sphere contains a straight line segment of length 2. In fact, we prove that the rectangular modulus attains its upper bound iff the unit sphere contains a straight line segment of length 2. 
 We further prove that if the dimension of the space $\mathbb{X}$ is finite then $\mu(\mathbb{X})$ is attained.  We also prove that a normed linear space is an inner product space iff sup$\{\frac{1+|t|}{\|y+tx\|}$: $x,y \in S_{\mathbb{X}}$ with $x\bot_By\} \leq \sqrt{2}$ $\forall t$ satisfying $|t|\in (3-2\sqrt{2},\sqrt{2}+1)$.\\

\section{Main Results}
 Let \[ \mbox{}~~\mu(x,y) = \frac{\|x\|+\|y\|}{\|x+y\|},  x, y  \in \mathbb{X} - \{0\}, x \perp_B y.  \]
Then \[  \mu(\mathbb{X}) = \sup \{ \mu(x,y) : x, y  \in \mathbb{X} - \{0\}, x \perp_B y\}. \]  
As the rectangular constant is defined by taking the supremum of $ \mu(x,y)  $  over the set $\{ (x,y) \in (\mathbb{X} - \{0\}) \times (\mathbb{X} - \{0\}) : x \perp_B y \} $ so it is natural to ask whether it is attained i.e., whether there exists some $ x_0, y_0  \in \mathbb{X} - \{0\}, x_0 \perp_B y_0 $ such that  
$ \mu(\mathbb{X}) = \mu(x_0,y_0).$ \\
The mapping 
\[ (x,y) \longrightarrow \mu(x,y) \] 
is continuous from $(\mathbb{X} - \{0\}) \times (\mathbb{X} - \{0\}) $ to $\mathbb{R}$ but the domain $\{ (x,y) \in (\mathbb{X} - \{0\}) \times (\mathbb{X} - \{0\}) : x \perp_B y \} $ is not  compact and so we can not trivially say that the rectangular constant is attained, even if the space is finite dimensional. We here show that the rectangular constant is attained in a finite dimensional normed linear space. We first prove a lemma which we will use repeatedly.
\begin{lemma}
Let $\mathbb{X}$ be a normed linear space. If $x, y \in S_{\mathbb{X}}$ such that $\|tx+(1-t)y\|=1 \forall ~t \in (0,1)$, then $x \bot_B (y-x)$.
\end{lemma}

\noindent\textbf{Proof.}
 For this it is sufficient to show that  $\|x + k (y-x)\| \geq 1 \forall k$ such that $|k|<1$.\\
If $0<k<1$ then  $\|x+k(y-x)\|= \|(1-k)x+ky\|=1$ and if  $-1<k<0$, then $\|x + k (y-x)\| = \|(1-k)x+ ky\| \geq \mid \|(1-k)x \| - \|ky\| \mid = 1. $

\begin{theorem}
Rectangular constant $\mu(\mathbb{X})$ is attained in any finite dimensional normed linear space $ \mathbb{X}. $
\end{theorem}

\noindent\textbf{Proof.}
 It is easy to check that the rectangular constant  $\mu(\mathbb{X}) $ can be written as 
\[ \mu(\mathbb{X}) = \sup \Big\{\frac{1+|t|}{\|y+tx\|} : x, y \in S_{\mathbb{X}}, x \bot_B y, t \in \mathbb{R}\Big\}.\]
Let $M = \{(x,y)\in S_{\mathbb{X}}\times S_{\mathbb{X}} : x \bot_B y \}$. Clearly $M$ is compact as $\mathbb{X}$ is finite dimensional. We consider the following two cases.\\
\noindent\textbf{Case 1.} Assume that  $\forall (x,y) \in M$, $\|sx+(1-s)y\|<$ $1$ for all $s \in (0,1)$.\\
 Then $1-\frac{\|x+y\|}{2}>$ $0$ $\forall(x,y) \in M$. Define $f: M\longrightarrow \mathbb{R} $ by
\[ f(x,y) = 1 - \frac{\|x+y\|}{2}. \]
Then $f$ is a continuous function on the compact set $M$ and so attains its infimum.\\
Let $\delta = \inf \{1-\frac{\|x+y\|}{2}: (x,y)\in M\}$, then $\delta> 0$. 
Now for  each $(x,y)\in M$, define $g_{xy}:\mathbb{R}\longrightarrow \mathbb{R}$ by 
\[ g_{xy}(t)=\frac{1+|t|}{\|y+tx\|}.\]
 Then $ 1 \leq  g_{xy}(t) \leq \frac{1+|t|}{|t|}$.\\ 
Now $g_{xy}(1) = \frac{2}{\|x+y\|} \geq \frac{1}{1-\delta} > 1$ for all $(x,y) \in M.$\\
Then we can find $k > 1$ such that $ g_{xy}(t) < g_{xy}(1) $ for all $|t| \geq k$. \\
This is true for all $(x, y) \in M$.\\
Define a function $g:M\times [-k,k] \longrightarrow \mathbb{R}$ 
\[ g(x,y,t) = \frac{1+|t|}{\|y+tx\|}\]
 which is continuous on the compact domain and  so attains its supremum.\\
\noindent\textbf{Case 2.}
Assume that $\exists (x,y) \in M$ such that $\|sx+(1-s)y\|=1$ for some $s$ in $(0,1)$. Then by the convexity of the unit sphere $\|sx+(1-s)y\| = 1$ $\forall s \in (0,1)$ as $\|x\| = \|y\| = 1$.\\
Then $x \bot_B (y-x)$ by the previous lemma.
Now $\mu(\mathbb{X}) \geq \frac{1+\|y-x\|}{\|x+(y-x)\|} = 1+\|y-x\|$.\\ Now $x\bot_B y$. 
If $\|y-x\| > 1$, then $\mu(\mathbb{X}) > 2 $.\\
Our claim is that $\mu(\mathbb{X})$ is then attained.\\
If $|t| > 3$ then $\|y+tx\| \geq |t|-1 > \frac{1+|t|}{2}$ and so  $\frac{1+|t|}{\|y+tx\|} < 2$.\\
 As $\mu(\mathbb{X}) > 2$ so we need to consider those $t$'s for which $ \mid t \mid \leq 3$.
Define a function $f:M\times [-3,3] \longrightarrow \mathbb{R}$ by
\[ f(x,y,t) = \frac{1+|t|}{\|y+tx\|}\] which is continuous on the compact domain and so attains its supremum. Hence $\mu(\mathbb{X})$ is attained.\\
If $\|y-x\| = 1$, then $\mu(\mathbb{X}) \geq 2$ and in case of $\mu(\mathbb{X}) = 2$, we see that \[\frac{\|y\|+\|-x\|}{\|y-x\|} = 2 \] and $\mu(\mathbb{X})$ is attained.\\
\noindent This completes the proof.
\begin{remark}
Although we could not give an example to show that the rectangular constant may not be attained in an infinite dimensional normed linear space, existence of such a space seems more probable to us.
\end{remark}
\section{On the characterisation of attainment of upper bound of rectangular constant}
\noindent We next prove that the rectangular constant $ \mu(\mathbb{X}) = 3 $ iff the unit sphere contains a straight line segment of length 2. In fact, we prove that the rectangular modulus attains its upper bound iff the unit sphere contains a straight line segment of length 2. For this we first prove the following theorem.
\begin{theorem}
Let $\mathbb{X}$ be a finite dimensional normed linear space. Let $l$ be a real number such that
\[ \|u-v\| < l ~\forall ~u,v \in S_\mathbb{X}  ~\mbox{with}~\|tu+(1-t)v\|=1 \forall ~t \in [0,1] \]
Then for any $x,y \in S_\mathbb{X}$ with $x\bot_By$ we have $\|x + \lambda y\|>1$ for all $\lambda$ with $|\lambda| \geq l.$
\end{theorem}
\noindent\textbf{Proof.} If possible let there be $x,y \in S_\mathbb{X}$ with $ x\bot_By $   such that $\|x+\lambda_0y\|= 1$ for some $\lambda_0$ with $|\lambda_0| \geq l$.\\
Now for  $t \in (0,1) $  we have
\begin{eqnarray*}
1 & \geq & t\|x\|+(1-t)\|x+\lambda_0y\| \\
  & \geq & \|tx+(1-t)(x+\lambda_0y)\|\\
	& = & \|x+\lambda_0(1-t)y\| \\
	&\geq & 1 
\end{eqnarray*}
$\Rightarrow \|tx+(1-t)(x+\lambda_0y)\|=1$ for all $t \in (0,1)$.\\
 So  $\|x -( x + \lambda_0y ) \|=\|\lambda_0y\|=|\lambda_0| \geq l$, this is a contradiction.\\
Hence  $\|x+\lambda y\|>1$ for all $\lambda$ with $|\lambda| \geq l$.\\
This completes the proof. \\
The converse part of the above theorem also holds as follows from the following theorem.
\begin{theorem}
Let $\mathbb{X}$ be a finite dimensional normed linear space. Let $l$ be a real number such that
for any $x,y \in S_\mathbb{X}$ with $x\bot_By$ we have $\|x + \lambda y\|>1$ for all $\lambda$ with $|\lambda| \geq l.$ Then 
\[ \|u-v\| < l ~\forall ~u,v \in S_\mathbb{X}  ~\mbox{with}~\|tu+(1-t)v\|=1 \forall ~t \in [0,1] \]
\end{theorem} 
\noindent\textbf{Proof.} 
If possible let there be  two elements $u,v \in S_\mathbb{X}$  satisfying $\|tu+(1-t)v\|=1$ $\forall t \in (0,1)$ such that  $\|u-v\| \geq l$.\\
Then by Lemma $2.1$ we have $u\bot_B(v-u)$ and so by homogenity of Birkhoff-James orthogonality $u\bot_B\frac{v-u}{\|v-u\|}.$ 
Let $x = u, y = \frac{v-u}{\|v-u\|}$ and $ \lambda = \| u - v\| .$\\
 Then $ \mid \lambda \mid \geq l $ and $\| x + \lambda y \|=\|u+v-u\| = 1, $ which is a contradiction to our assumption.\\
This completes the proof.\\
We next prove the theorem which shows that the $*$-rectangular modulus attains its upper bound iff the unit sphere contains a straight line segment of length 2.
\begin{theorem}
Let $\mathbb{X}$ be a finite dimensional normed linear space. Then $*$-rectangular modulus $\mu^*_{\mathbb{X}}(\lambda) = \lambda+2 $ iff $\exists$ two points $u',v' \in S_{\mathbb{X}} $ such that $\|tu'+(1-t)v'\|=1$ for all $t \in (0,1)$ and $\|u'-v'\|=2$.
\end{theorem}
\noindent\textbf{Proof.} We prove the necessary part by the method of contradiction i.e., we assume $\mu^*_{\mathbb{X}}(\lambda) = \lambda+2 $ and 
$ \| u' - v'\| < 2 $ whenever  $ u',v' \in S_{\mathbb{X}}$ satisfies $\|tu'+(1-t)v'\|=1$ for all $t \in(0,1).$ \\
It is easy to check that the $*$-rectangular modulus $ \mu^*_{\mathbb{X}}(\lambda)$ can be written as 
 \[ \mu^*_{\mathbb{X}}(\lambda) = \sup \Big \{\frac{\lambda+t}{\|u+tv\|} : u\bot_Bv,u,v \in S_{\mathbb{X}}, t>0 \Big \} .\]
If $ 0 < t \leq 1$ then $\frac{\lambda+t}{\|u+tv\|} \leq 1+\lambda$ and if  $1<t<2$ then $\frac{\lambda+t}{\|u+tv\|}<\lambda+2$.\\
If $2 \leq t <2+\frac{1}{\lambda}$, then $\|u+tv\|\geq t-1$ and $\|u+tv\|>1$  by the  Theorem 3.1. \\
So $\|u+tv\|>1.\frac{\lambda}{\lambda+2}+(t-1)\frac{2}{\lambda+2} \geq \frac{\lambda+t}{\lambda+2}$, which implies\\
$\frac{\lambda+t}{\|u+tv\|}< \lambda+2$.\\
Finally let $t \geq 2+\frac{1}{\lambda}$.\\
Then $\|u+tv\| \geq t-1 \geq \frac{\lambda+t}{\lambda+1}$, which implies $\frac{\lambda+t}{\|u+tv\|} \leq \lambda+1$. \\
So we see that sup$\{\frac{\lambda+t}{\|u+tv\|}:u\bot_Bv,u,v \in S_{\mathbb{X}}, 0<t \leq 1 \} \leq 1+\lambda$ and\\
sup$\{\frac{\lambda+t}{\|u+tv\|}:u\bot_Bv, u,v \in S_{\mathbb{X}}, t\geq 2+\frac{1}{\lambda}\} \leq 1+\lambda$.\\
Let $M=\{(u,v)\in S_X\times S_{\mathbb{X}}:u\bot_Bv \}$ which is compact, since ${\mathbb{X}}$ is finite dimensional.\\
Define a function $f:M\times[1,2+\frac{1}{\lambda}]\longrightarrow \mathbb{R}$ by $f(u,v,t)=\frac{\lambda+t}{\|u+tv\|}$ which is continuous on a compact set, so attains its supremum. So sup$\{\frac{\lambda+t}{\|u+tv\|}:u\bot_Bv,u,v \in S_{\mathbb{X}},t>0\} < 2+\lambda$.\\
 
\noindent Conversely, suppose that there exist two points $u,v \in S_{\mathbb{X}}$ such that $\|tu+(1-t)v\|=1$ for all $t \in (0,1)$ and $\|u-v\|=2$. Now by Lemma $2.1$ $u\bot_B(v-u)$ and so $u\bot_B\frac{v-u}{\|v-u\|}$, as previously.\\
Now, $\frac{\lambda+\|v-u\|}{\|u+\|v-u\|\frac{v-u}{\|v-u\|}\|}=\lambda+\|v-u\|=\lambda+2$.\\
Hence $\mu^*_\mathbb{X}(\lambda) =$ sup$\{\frac{\lambda+t}{\|u+tv\|}:u\bot_Bv,u,v\in S_{\mathbb{X}}\} = 2+\lambda$.\\
This completes the proof.\\

\begin{cor}Let $\mathbb{X}$ be a finite dimensional normed linear space. Then rectangular modulus $\mu_\mathbb{X}(\lambda) =$ max$\{2+\lambda,1+2\lambda\}$ iff $\exists$ two points $u,v \in S_{\mathbb{X}}$ such that $\|tu+(1-t)v\|=1$ for all $t \in (0,1)$ and $\|u-v\|=2$.\\
\end{cor}

\noindent\textbf{Proof.}
Since $\mu_{\mathbb{X}}(\lambda) =$ max$\{\mu^*_{\mathbb{X}}(\lambda),\lambda \mu^*_{\mathbb{X}}(\frac{1}{\lambda})\}$, the result follows from the last theorem and the fact that $\mu^*_{\mathbb{X}}(\lambda) = \lambda+2$ iff $\lambda\mu^*_{\mathbb{X}}(\frac{1}{\lambda}) = 1+2\lambda$.

\begin{cor}
Let $\mathbb{X}$ be a finite dimensional normed linear space. Then rectangular constant $\mu(\mathbb{X}) = 3$ iff $\exists$ two points $u,v \in S_{\mathbb{X}}$ such that $\|tu+(1-t)v\|=1$ for all $t \in (0,1)$ and $\|u-v\|=2$.
\end{cor}

\noindent\textbf{Geometric interpretation.} The above result shows that the rectangular constant of the space is equal to $3$ iff the unit sphere of space contains a straight line segment of length $2.$ \\

\noindent We next prove that if $\mathbb{X}$ is a two dimensional polyhedral space ( i.e., a space with finitely many extreme points ) then in finding $\mu({\mathbb{X}}) $ it is sufficient to consider the supremum over all orthogonal  $x, y \in \mathbb{X} - \{0\} $ where $x$ is an extreme point of $S_{\mathbb{X}}.$
\begin{theorem}
Let $\mathbb{X}$ be a two-dimensional real polyhedral space. Then for any $x, y \in \mathbb{X}-\{0\} $ with $x\bot_By,$ 
$\mu(x,y)= \frac{\|x\|+\|y\|}{\|x+y\|} \leq \frac{\|u\|+\|v\|}{\|u+v\|}$ where $u \in S_{\mathbb{X}}$ is an extreme point and $u\bot_Bv$.
\end{theorem}
\noindent\textbf{Proof.} We first prove the following two lemmas: 
\begin{lemma}
Let $\mathbb{X}$ be a finite dimensional polyhedral space. Let $x_{1}, x_{2}\in S_{\mathbb{X}}$ are two adjacent extreme points and let $x= tx_{1}+(1-t)x_{2}$,
for some $t\in (0,1)$. Then $x\bot_By \Rightarrow x_{1}\bot_By$ and $x_{2}\bot_By$.
\end{lemma}

\noindent\textbf{Proof.} Assume $ x \bot_B y.$ 
By Theorem 2.1 of James \cite{6} we can find a linear functional $f$ on ${\mathbb{X}}$ such that $\|f\| = f(x) = 1 $  and $ f(y) = 0.$\\
Now $|f(x_{1})|\leq \|x_{1}\|= 1$ and $|f(x_{2})|\leq \|x_{2}\|= 1$.\\
If $|f(x_{1})|< 1$ then
\begin{eqnarray*}
1 & = & \|tx_{1}+(1-t)x_{2}\|\\
 & = & t\|x_{1}\|+(1-t)\|x_{2}\|\\
 & > & t|f(x_{1})|+(1-t)|f(x_{2})|\\
 & \geq &  |f(tx_{1}+(1-t)x_{2})|\\
  & = &|f(x)|,~~ \mbox{which is a contradiction.}
\end{eqnarray*} 
So $|f(x_{1})|= 1$. Similarly $|f(x_2)|= 1$. As $f(y) = 0 $ so by the same Theorem 2.1 of James \cite{6} we get $x_{1}\bot_By$ and $x_{2}\bot_By$. \\
Thus $ x \bot_B y \Rightarrow x_{1}\bot_By$ and $x_{2}\bot_By$. This completes the proof of first lemma.

\begin{lemma}
Let $x, y \in {\mathbb{X}}$ be such that $x\bot_By$. Then $\|x+\lambda_{1} y\| \geq \|x+\lambda_{2} y\|$, if $\lambda_{1} > \lambda_{2} > 0$ or $\lambda_{1} < \lambda_{2} < 0$.
\end{lemma}
\noindent\textbf{Proof.}
Let $\lambda_{1} > \lambda_{2} > 0$, then we can write $\lambda_{1}= \lambda_{2}+c$, for some $c > 0$. To show \[ \|x+\lambda_{1}y\|= \|x+\lambda_{2}y+cy\| \geq \|x+\lambda_{2}y\|\]\\
If $\|x+\lambda_{2}y\|= \|x\|$, then \[\|x+\lambda_{1}y\| \geq \|x\|= \|x+\lambda_{2}y\|\] and we are done.\\
Let $\|x+\lambda_{2}y\| > \|x\|$.\\
We know, for ($x+\lambda_{2}y$) and $y$ either $\|x+\lambda_{2}y+\lambda y\| \geq \|x+\lambda_{2}y\| \forall \lambda \geq 0$ or $\|x+\lambda_{2}y+\lambda y\| \geq \|x+\lambda_{2}y\| \forall \lambda \leq 0$.\\
Put $\lambda= -\lambda_{2} < 0$, then $\|x\|= \|x+\lambda_{2}y-\lambda_{2}y\| < \|x+\lambda_{2}y\|$\\
So $\|x+\lambda_{2}y+cy\| \geq \|x+\lambda_{2}y\|$ i.e. $\|x+\lambda_{1}y\| \geq \|x+\lambda_{2}y\|$.\\
Similarly we can prove the result for $\lambda_{1} < \lambda_{2} < 0$. This completes the proof of the lemma.\\

\noindent\textbf{Proof of the theorem continued:}\\
Let $x', y' \in \mathbb{X}- \{0\}$ and $x' \bot_B y'$. Then $\mu(x', y')= \frac{\|x'\|+\|y'\|}{\|x'+y'\|}= \frac{1+|\lambda|}{\|x+\lambda y\|}$, for some $x,y \in S_{{\mathbb{X}}}$ with $x\bot_By$. \\
If $x\in S_{{\mathbb{X}}}$ is an extreme point, then we are done.\\
Suppose $x\in S_{{\mathbb{X}}}$ is not an extreme point. Then $x= tx_{1}+(1-t)x_{2}$, for some $t\in(0,1)$ and for two adjacent extreme points $x_{1}, x_{2}\in S_{{\mathbb{X}}}$.\\
Let $|\lambda| \leq \|x_{1}-x_{2}\| $\\
By Lemma $2.1$ $x_{1}\bot_B(x_{2}-x_{1})$. Now
 \[\frac{1+|\lambda|}{\|x+\lambda y\|}\leq 1+|\lambda|\leq 1+\|x_{2}-x_{1}\|\leq \frac{\|x_{1}\|+\|x_{2}-x_{1}\|}{\|x_{1}+x_{2}-x_{1}\|}\].

\noindent Let $|\lambda|> \|x_{1}-x_{2}\|$.
As $x\bot_By$ so $x_{1}\bot_By$ and $x_{2}\bot_By$ by  Lemma $3.7$. Then $x_{1}= x+\lambda_{1} y$ and $x_{2}= x+\lambda_{2} y$ for some $\lambda_{1}, \lambda_{2} \in \mathbb{R} $.\\
Claim: $\lambda_{1}, \lambda_{2}$ are of opposite sign and $|\lambda_{1}|, |\lambda_{2}|< \|x_{1}-x_{2}\|$.\\
If possible, let $\lambda_{1}> 0, \lambda_{2}> 0$. Then 
\begin{eqnarray*}
x_{1} = x+\lambda_{1} y & = & tx_{1}+(1-t)x_{2}+ \lambda_{1}y \\
\Rightarrow (1-t)x_{1} &=  & (1-t)x_{2}+ \lambda_{1} y\\ 
\Rightarrow x_{1}-x_{2} & = & \frac{\lambda_{1}y}{1-t}\ldots ~~~~~~~~~~~~~~~~~~~~~~~~~(i)\\ 
\Rightarrow \|x_{1}-x_{2}\| & = & \frac{|\lambda_{1}|}{1-t} > |\lambda_{1}|
\end{eqnarray*}
Proceeding similarly from  $x_{2}= x+\lambda_{2} y $,\\
 we get $x_{2}-x_{1}= \frac{\lambda_{2} y}{t} \ldots ~~~~~~~~~~~~~~~~~~~~~~~~~(ii)\\
\Rightarrow |\lambda_{2}|< \|x_{1}-x_{2}\|$\\
From $(i)$ and $(ii)$ we get $\frac{\lambda_{1}}{t-1} = \frac{\lambda_{2}}{t}\\
                             \Rightarrow t= \frac{\lambda_{2}}{\lambda_{2}-\lambda_{1}}$\\
If $\lambda_{1} > \lambda_{2}$, then $t < 0$ and if $\lambda_{2} > \lambda_{1}$, then $t > 1$, contradiction in both cases. Similar case for $\lambda_{1} < 0$, $\lambda_{2} < 0$. \\

\noindent Let us first assume that  $\lambda_1 > 0$, $\lambda_{2} < 0$ and $\lambda > 0$.\\
As $|\lambda_{2}| < \|x_{1}-x_{2}\| < |\lambda|$, we have $0 < \lambda+\lambda_{2} < \lambda$.\\
Then by Lemma $3.8$, $\|x_{2}+\lambda y\|= \|x+\lambda_{2}y+\lambda y\| \leq \|x+\lambda y\|$.\\
Let us next assume that $\lambda_1 > 0$, $\lambda_{2} < 0$ and $\lambda < 0$, then $\lambda < \lambda+\lambda_{1} < 0$.\\
So $\|x_{1}+\lambda y\|= \|x+\lambda_{1}y+\lambda y\|= \|x+(\lambda_{1}+\lambda)y\| \leq \|x+\lambda y\|$.\\
Similarly, if $\lambda_{1} < 0$, $\lambda_{2} > 0$, then,\\
$\|x_{1}+\lambda y\| \leq \|x+\lambda y\|$, when $\lambda > 0 \\
\|x_{2}+\lambda y\| \leq \|x+\lambda y\|$, when $\lambda < 0$.\\
Now $\mu(x', y')= \frac{1+|\lambda|}{\|x+\lambda y\|} \leq \frac{1+|\lambda|}{\|x_{i}+\lambda y\|}$, $i \in \{1,2\} \\  
\Rightarrow \mu(x', y') \leq \frac{\|x_{i}\|+\|\lambda y\|}{\|x_{i}+\lambda y\|}$, $i \in \{1,2\}$.\\
This completes the proof.	\\

\noindent We next give an example showing that $v$ in Theorem $3.6$ can not be taken as either an extreme point or a scalar multiple of an extreme point on the unit sphere.
\begin{example}
 Consider the normed linear space $(\mathbb{R}^2,\|~~\|_\infty)$, which is a two-dimensional polyhedral space. The only extreme points of the unit sphere are $\pm(1,1), \pm(-1,1)$. We first note that $(1,1)\bot_B k(1,-1)$ for any real k and $(1,1)\bot_B(-2,0)$. We also see that $\frac{\|(1,1)\|+\|k(1,-1)\|}{\|(1,1)+k(1,-1)\|} = 1 < \frac{\|(1,1)\|+\|(-2,0)\|}{\|(1,1)+(-2,0)\|} = 3$. Similarly considering other extreme points of the unit sphere and their scalar multiples, we conclude that $v$ in the previous theorem can not be taken as an extreme point in this case.
\end{example}   
\noindent We next prove the theorem which characterises an inner product space. Joly \cite{8} and del Rio and Benitez\cite{3} proved that a finite dimensional normed linear space is an inner product space iff $ sup \{ \frac{1 + \mid \lambda \mid }{\| y + \lambda x\|} : x \bot_B y, x, y \in S_{\mathbb{X}}, \lambda \in \mathbb{R} \} = \sqrt{2}. $  Here we prove that 
a finite dimensional normed linear space is an inner product space iff $ sup \{ \frac{1 + \mid \lambda \mid }{\| y + \lambda x\|} : x \bot_B y, x, y \in S_{\mathbb{X}}, \mid \lambda \mid \in (3 - 2\sqrt{2}, \sqrt{2} + 1) \} \leq \sqrt{2}. $
\begin{theorem}
A finite dimensional normed linear space $\mathbb{X}$ is an inner-product space iff $x,y \in S_{\mathbb{X}}$ with $x\bot_By$ we have $\frac{1+|\lambda|}{\|y+\lambda x\|} \leq \sqrt{2}$ $\forall \lambda$ satisfying $|\lambda|\in (3-2\sqrt{2},\sqrt{2}+1)$.
\end{theorem}
\noindent\textbf{Proof.} The necessary part of the theorem follows easily. For the sufficient part we proceed as follows. \\
As $x\bot_By$, we have $\|y+\lambda x\| \geq |\lambda| .\\
\Rightarrow \frac{1+|\lambda|}{\|y+\lambda x\|} \leq \frac{1+|\lambda|}{|\lambda|} \leq \sqrt{2}$ $\forall \lambda$ satisfying $|\lambda| \geq \sqrt{2}+1$.\\
Again, $\|y+\lambda x\| \geq |1-|\lambda||$\\
Let $\|y+\lambda x\| \geq 1-|\lambda|\\
\Rightarrow \frac{1+|\lambda|}{\|y+\lambda x\|} \leq \frac{1+|\lambda|}{1-|\lambda|} \leq \sqrt{2}$ $\forall \lambda$ satisfying $|\lambda| \leq 3-2\sqrt{2}$.\\
Let $\|y+\lambda x\| \geq |\lambda|-1\\
\Rightarrow \frac{1+|\lambda|}{\|y+\lambda x\|} \leq \frac{1+|\lambda|}{|\lambda|-1} \leq \sqrt{2}$ for all $\lambda$ satisfying $|\lambda| \geq 3+2\sqrt{2} > \sqrt{2}+1$.\\
So if $\frac{1+|\lambda|}{\|y+\lambda x\|} \leq \sqrt{2}$ for all $\lambda$ with $|\lambda| \in (3-2\sqrt{2},\sqrt{2}+1)$, then $\frac{1+|\lambda|}{\|y+\lambda x\|} \leq \sqrt{2} \forall \lambda \in \mathbb{R}$. So $\mathbb{X}$ is an inner product space.\\
This completes the proof.

\end{document}